# q-Chebyshev polynomials

*Johann Cigler*

Fakultät für Mathematik, Universität Wien

*johann.cigler@univie.ac.at*

**Abstract**

In this overview paper a direct approach to $q-$Chebyshev polynomials and their elementary properties is given. Special emphasis is placed on analogies with the classical case. There are also some connections with $q-$tangent and $q-$Genocchi numbers.

## 0. Introduction

Waleed A. Al Salam and Mourad E.H. Ismail [1] found a class of polynomials which can be interpreted as $q-$analogues of the bivariate Chebyshev polynomials of the second kind. These are essentially the polynomials $U_n(x,s,q)$ which will be introduced in (2.12). In [11] I also considered corresponding $q-$Chebyshev polynomials $T_n(x,s,q)$ of the first kind which will be defined in (2.6). Together these polynomials satisfy many $q-$analogues of well-known identities for the classical Chebyshev polynomials $T_n(x) = T_n(x,-1,1)$ and $U_n(x) = U_n(x,-1,1)$. For some of them it is essential that our polynomials depend on two independent parameters. This is especially true for (2.36) which generalizes the defining property $\left(x+\sqrt{x^2-1}\right)^n = T_n(x) + U_{n-1}(x)\sqrt{x^2-1}$ of the classical Chebyshev polynomials. Another approach to univariate $q-$analogues of Chebyshev polynomials has been proposed by Natig Atakishiyev et al. in [2], (5.3) and (5.4). In our terminology they considered the monic versions of the polynomials $T_n\left(x, -\frac{1}{\sqrt{q}}, q\right)$ and $U_n\left(x, -\sqrt{q}, q\right)$. Since $U_n(x, s^2, q) = s^n U_n\left(\frac{x}{s}, 1, q\right)$ and $T_n(x, s^2, q) = s^n T_n\left(\frac{x}{s}, 1, q\right)$ their definition also leads to the same bivariate polynomials $T_n(x,s,q)$ and $U_n(x,s,q)$.
Without recognizing them as $q-$analogues of Chebyshev polynomials some of these polynomials also appeared in the course of computing Hankel determinants as in [7] and [13].

The purpose of this paper is to give a direct approach to these polynomials and their simplest properties.

# 1. Some well-known facts about the classical Chebyshev polynomials

Let me first state some well-known facts about those aspects of the classical Chebyshev polynomials (cf. e.g. [15]) and their bivariate versions for which we will give $q-$analogues.

The (classical) *Chebyshev polynomials of the first kind* $T_n(x)$ satisfy the recurrence

$$T_n(x) = 2xT_{n-1}(x) - T_{n-2}(x) \tag{1.1}$$

with initial values $T_0(x) = 1$ and $T_1(x) = x.$
For $x = 1$ this reduces to

$$T_n(1) = 1. \tag{1.2}$$

The (classical) *Chebyshev polynomials of the second kind* $U_n(x)$ satisfy the same recurrence

$$U_n(x) = 2xU_{n-1}(x) - U_{n-2}(x) \tag{1.3}$$

but with initial values $U_{-1}(x) = 0$ and $U_0(x) = 1,$ which gives $U_1(x) = 2x.$
As special values we note that
$$U_n(1) = n+1. \tag{1.4}$$

These polynomials are related by the identity

$$\left(x + \sqrt{x^2-1}\right)^n = T_n(x) + U_{n-1}(x)\sqrt{x^2-1}, \tag{1.5}$$

which in turn implies
$$T_n(x)^2 - (x^2-1)U_{n-1}(x)^2 = 1. \tag{1.6}$$

**Remark 1.1**

For $x = \cos\vartheta$ identity (1.5) becomes
$$\cos n\vartheta + i\sin n\vartheta = \left(\cos\vartheta + i\sin\vartheta\right)^n = T_n\left(\cos\vartheta\right) + iU_{n-1}(\cos\vartheta)\sin\vartheta$$
or equivalently

$$T_n(\cos\vartheta) = \cos n\vartheta$$
$$U_n(\cos\vartheta) = \frac{\sin(n+1)\vartheta}{\sin\vartheta}. \tag{1.7}$$

This is the usual approach to the Chebyshev polynomials. Identity (1.6) reduces to

$$\cos^2 n\vartheta + \sin^2 n\vartheta = 1. \tag{1.8}$$

Unfortunately it seems that this aspect of the Chebyshev polynomials has no simple $q-$analogue.

The Chebyshev polynomials are *orthogonal polynomials*. As is well-known (cf. e.g. [4]) a sequence $(p_n(x))_{n \geq 0}$ of polynomials with $p_0(x) = 1$ and $\deg p_n = n$ is called *orthogonal with respect to a linear functional* $\Lambda$ on the vector space of polynomials if $\Lambda(p_m p_n) = 0$ for $m \neq n$. The linear functional is uniquely determined by $\Lambda(p_n) = [n = 0]$. Here $[P]$ denotes the Iverson symbol defined by $[P] = 1$ if property $P$ is true and $[P] = 0$ otherwise. The values $\Lambda(x^n)$ are called moments of $\Lambda$.

Let $P_n(x)$ denote the monic polynomials corresponding to $p_n(x)$ and $a(n,k)$ be the uniquely determined coefficients in

$$\sum_{k=0}^{n} a(n,k) P_k(x) = x^n. \tag{1.9}$$

Then $a(n,0) = \Lambda(x^n)$ and more generally $a(n,k) = \dfrac{\Lambda(x^n P_k(x))}{\Lambda(P_k(x)^2)}$.

By Favard's theorem **there exist numbers $s(n), t(n)$ such that the three-term recurrence**

$$P_n(x) = (x - s(n-1)) P_{n-1}(x) - t(n-2) P_{n-2}(x) \tag{1.10}$$

holds.

Therefore the coefficients $a(n,k)$ satisfy

$$\begin{aligned} a(0,j) &= [j = 0] \\ a(n,0) &= s(0) a(n-1,0) + t(0) a(n-1,1) \\ a(n,j) &= a(n-1, j-1) + s(j) a(n-1, j) + t(j) a(n-1, j+1). \end{aligned} \tag{1.11}$$

This can be used to compute the moments $a(n,0) = \Lambda(x^n)$.
If the moments are known then the corresponding orthogonal polynomials $P_n(x)$ are given by

$$P_n(x) = \frac{1}{\det(\Lambda(x^{i+j}))_{i,j=0}^{n-1}} \det \begin{pmatrix} \Lambda(x^0) & \Lambda(x^1) & \cdots & \Lambda(x^{n-1}) & 1 \\ \Lambda(x^1) & \Lambda(x^2) & \cdots & \Lambda(x^n) & x \\ \Lambda(x^2) & \Lambda(x^3) & \cdots & \Lambda(x^{n+1}) & x^2 \\ \vdots & & & & \vdots \\ \Lambda(x^n) & \Lambda(x^{n+1}) & \cdots & \Lambda(x^{2n-1}) & x^n \end{pmatrix}. \tag{1.12}$$

So the knowledge of the polynomials $P_n(x)$ is equivalent with the knowledge of $s(n)$ and $t(n)$ and this is again equivalent with the knowledge of the moments.

Since $\frac{1}{\pi}\int_{-1}^{1}\frac{T_n(x)}{\sqrt{1-x^2}}dx = \frac{1}{\pi}\int_0^\pi \cos(n\vartheta)d\vartheta = [n=0]$ for the polynomials $T_n(x)$ the corresponding linear functional $L$ is given by the integral

$$L(p(x)) = \frac{1}{\pi}\int_{-1}^{1}\frac{p(x)}{\sqrt{1-x^2}}dx \tag{1.13}$$

and

$$L(T_n^2) = \frac{1}{\pi}\int_{-1}^{1}\frac{T_n(x)^2}{\sqrt{1-x^2}}dx = \begin{cases} 1 & \text{if } n=0 \\ \frac{1}{2} & \text{if } n>0 \end{cases} \tag{1.14}$$

The corresponding moments are

$$L(x^{2n}) = \frac{1}{\pi}\int_{-1}^{1}\frac{x^{2n}}{\sqrt{1-x^2}} = \frac{1}{2^{2n}}\binom{2n}{n} \tag{1.15}$$

and $L(x^{2n+1}) = 0$.

For the polynomials $U_n(x)$ we get from

$$\frac{2}{\pi}\int_{-1}^{1}U_n(x)\sqrt{1-x^2}dx = \frac{2}{\pi}\int_0^\pi \sin((n+1)\vartheta)\sin\vartheta d\vartheta = [n=0]$$

that the corresponding linear functional $M$ satisfies

$$M(p_n) = \frac{2}{\pi}\int_{-1}^{1}p(x)\sqrt{1-x^2}dx \tag{1.16}$$

and

$$M(U_n^2) = \frac{2}{\pi}\int_{-1}^{1}U_n(x)^2\sqrt{1-x^2}dx = 1. \tag{1.17}$$

The corresponding moments are

$$M(x^{2n}) = \frac{2}{\pi}\int_{-1}^{1}x^{2n}\sqrt{1-x^2}dx = \frac{1}{2^{2n}}\frac{1}{n+1}\binom{2n}{n} \tag{1.18}$$

and $M(x^{2n+1}) = 0$.

As already mentioned in the introduction for our $q-$ analogues we need bivariate Chebyshev polynomials.

The *bivariate Chebyshev polynomials* $T_n(x,s)$ *of the first kind* satisfy the recurrence

$$T_n(x,s) = 2xT_{n-1}(x,s) + sT_{n-2}(x,s) \tag{1.19}$$

with initial values $T_0(x,s) = 1$ and $T_1(x,s) = x$.
Of course $T_n(x) = T_n(x,-1)$.

They have the determinant representation

$$T_n(x,s) = \det\begin{pmatrix} x & s & 0 & \cdots & 0 & 0 \\ -1 & 2x & s & \cdots & 0 & 0 \\ 0 & -1 & 2x & \cdots & 0 & 0 \\ \vdots & \vdots & \vdots & \vdots & \vdots & \vdots \\ 0 & 0 & 0 & \cdots & 2x & s \\ 0 & 0 & 0 & \cdots & -1 & 2x \end{pmatrix}. \tag{1.20}$$

The *bivariate Chebyshev polynomials of the second kind* $U_n(x,s)$ satisfy the same recurrence

$$U_n(x,s) = 2xU_{n-1}(x,s) + sU_{n-2}(x,s) \tag{1.21}$$

but with initial values $U_0(x,s) = 1$ and $U_1(x,s) = 2x$.

Their determinant representation is

$$U_n(x,s) = \det\begin{pmatrix} 2x & s & 0 & \cdots & 0 & 0 \\ -1 & 2x & s & \cdots & 0 & 0 \\ 0 & -1 & 2x & \cdots & 0 & 0 \\ \vdots & \vdots & \vdots & \vdots & \vdots & \vdots \\ 0 & 0 & 0 & \cdots & 2x & s \\ 0 & 0 & 0 & \cdots & -1 & 2x \end{pmatrix}. \tag{1.22}$$

These polynomials are connected via

$$\left(x + \sqrt{x^2 + s}\right)^n = T_n(x,s) + U_{n-1}(x,s)\sqrt{x^2 + s}. \tag{1.23}$$

This also implies

$$T_n(x,s)^2 - (x^2 + s)U_{n-1}(x,s)^2 = (-s)^n. \tag{1.24}$$

The Chebyshev polynomials are intimately related with Fibonacci and Lucas polynomials

$$F_{n+1}(x,s) = \sum_{k=0}^{\lfloor \frac{n}{2} \rfloor} \binom{n-k}{k} s^k x^{n-2k} \tag{1.25}$$

and

$$L_n(x,s) = F_{n+1}(x,s) + sF_{n-1}(x,s) = \sum_{k=0}^{\lfloor \frac{n}{2} \rfloor} \frac{n}{n-k}\binom{n-k}{k}s^k x^{n-2k} \qquad (1.26)$$

for $n > 0$ (cf. e.g. [10]). Here as usual $L_0(x,s) = 2$.

More precisely the monic polynomials $T_0(x,s) = 1$ and $\dfrac{T_n(x,s)}{2^{n-1}}$ for $n > 0$ coincide with the modified Lucas polynomials

$$\frac{L_n^*(2x,s)}{2^n} = L_n^*\left(x, \frac{s}{4}\right). \qquad (1.27)$$

They are defined by $L_n^*(x,s) = L_n(x,s)$ for $n > 0$ and $L_0^*(x,s) = 1$ and satisfy a three-term recurrence with $s(n) = 0$, $t(0) = \dfrac{s}{2}$ and $t(n) = \dfrac{s}{4}$ for $n > 0$.

The moments can be obtained from the formula

$$\sum_{k=0}^{\lfloor \frac{n}{2} \rfloor} \binom{n}{k}(-s)^k L_{n-2k}^*(x,s) = x^n. \qquad (1.28)$$

The monic polynomials $\dfrac{U_n(x,s)}{2^n}$ are Fibonacci polynomials

$$\frac{U_n(x,s)}{2^n} = \frac{F_{n+1}(2x,s)}{2^n} = F_{n+1}\left(x, \frac{s}{4}\right). \qquad (1.29)$$

In this case the corresponding numbers $s(n)$ and $t(n)$ are $s(n) = 0$ and $t(n) = \dfrac{s}{4}$.

Here the moments can be obtained from

$$\sum_{k=0}^{\lfloor \frac{n+1}{2} \rfloor}\left(\binom{n}{k} - \binom{n}{k-1}\right)(-s)^k F_{n+1-2k}(x,s) = x^n. \qquad (1.30)$$

We shall also give $q-$ analogues of the following identities which express Chebyshev polynomials of odd order in terms of Chebyshev polynomials of even order:

$$T_{2n+1}(x) = \sum_{k=0}^{n}\binom{2n+1}{2k}(-1)^{n-k} t_{2n-2k+1} x^{2n+1-2k} T_{2k}(x) \qquad (1.31)$$

and

$$U_{2n+1}(x) = \sum_{k=0}^{n} \binom{2n+2}{2k} \frac{1}{2k+1} (-1)^{n-k} G_{2n-2k+2}(2x)^{2n-2k} U_{2k}(x). \tag{1.32}$$

Here the tangent numbers $(t_{2n+1})_{n\geq 0} = (1, 2, 16, 272, 7936, \cdots)$ and the Genocchi numbers $(G_{2n})_{n\geq 0} = (0, 1, 1, 3, 17, 155, 2073, \cdots)$ are given by their generating functions

$$\frac{e^z - e^{-z}}{e^z + e^{-z}} = \sum_{n \geq 0} (-1)^n \frac{t_{2n+1}}{(2n+1)!} z^{2n+1} \tag{1.33}$$

and

$$z \frac{e^z - e^{-z}}{e^z + e^{-z}} = \sum_{n \geq 0} (-1)^{n-1} 2^{2n-1} \frac{G_{2n}}{(2n)!} z^{2n}. \tag{1.34}$$

Note that

$$t_{2n+1} = \frac{2^{2n} G_{2n+2}}{n+1}. \tag{1.35}$$

## 2. q-analogues

We assume that $q \neq -1$ is a real number. All $q-$ identities in this paper reduce to known identities when $q$ tends to 1. We assume that the reader is familiar with the most elementary notions of $q-$ analysis (cf. e.g. [5]). The $q-$ binomial coefficients
$\begin{bmatrix} n \\ k \end{bmatrix} = \frac{[1][2]\cdots[n]}{[1]\cdots[k]\cdot[1]\cdots[n-k]}$ with $[n] = 1 + q + \cdots + q^{n-1}$ satisfy the recurrences

$$\begin{bmatrix} n \\ k \end{bmatrix} = q^k \begin{bmatrix} n-1 \\ k \end{bmatrix} + \begin{bmatrix} n-1 \\ k-1 \end{bmatrix} = \begin{bmatrix} n-1 \\ k \end{bmatrix} + q^{n-k} \begin{bmatrix} n-1 \\ k-1 \end{bmatrix}. \tag{2.1}$$

If we want to stress the dependence on $q$ we write $[n]_q$ and $\begin{bmatrix} n \\ k \end{bmatrix}_q$ respectively.

We also need the $q-$ Pochhammer symbol $(x; q)_n = (1-x)(1-qx)\cdots(1-q^{n-1}x)$ and the $q-$ binomial theorem in the form

$$(x; q)_n = \sum_{k=0}^{n} (-1)^k q^{\binom{k}{2}} \begin{bmatrix} n \\ k \end{bmatrix} x^k \tag{2.2}$$

or equivalently

$$p_n(x, y) = (x+y)(qx+y)\cdots(q^{n-1}x+y) = \sum_{k=0}^{n} q^{\binom{k}{2}} \begin{bmatrix} n \\ k \end{bmatrix} x^k y^{n-k}. \tag{2.3}$$

We denote by $e(z) = e(z,q) = \sum_{n\geq 0} \dfrac{z^n}{[n]!}$ the $q-$exponential function. It satisfies

$$\frac{1}{e(-z)} = \sum_{n\geq 0} q^{\binom{n}{2}} \frac{z^n}{[n]!}.$$

Since the Chebyshev polynomials are special cases of Fibonacci and Lucas polynomials it would be tempting to look for $q-$analogues related to the simplest $q-$analogues of Fibonacci and Lucas polynomials (cf. e.g. [10])

$$F_{n+1}(x,s,q) = \sum_{k=0}^{\lfloor \frac{n}{2} \rfloor} q^{k^2} \begin{bmatrix} n-k \\ k \end{bmatrix} s^k x^{n-2k},$$

$$L_n(x,s,q) = F_{n+1}(x,s,q) + sF_{n-1}(x,qs,q) = \sum_{k=0}^{\lfloor \frac{n}{2} \rfloor} q^{k^2-k} \frac{[n]}{[n-k]} \begin{bmatrix} n-k \\ k \end{bmatrix} s^k x^{n-2k},$$

$$Fib_{n+1}(x,s,q) = \sum_{k=0}^{\lfloor \frac{n}{2} \rfloor} q^{\binom{k+1}{2}} \begin{bmatrix} n-k \\ k \end{bmatrix} s^k x^{n-2k} \quad \text{and}$$

$$Luc_n(x,s,q) = Fib_{n+1}(x,s,q) + sFib_{n-1}(x,s,q) = \sum_{k=0}^{\lfloor \frac{n}{2} \rfloor} q^{\binom{k}{2}} \frac{[n]}{[n-k]} \begin{bmatrix} n-k \\ k \end{bmatrix} s^k x^{n-2k}.$$

But here we have no success. Though the polynomials $F_{n+1}(x,s,q)$ are orthogonal there are no closed forms for their moments. None of the other classes of polynomials satisfies a 3-term recurrence. So they cannot be orthogonal.

But it is interesting that for $Fib_{n+1}(x,s,q)$ and $Luc_n(x,s,q)$ the following analogues of (1.28) and (1.30)

$$\sum_{k=0}^{\lfloor \frac{n}{2} \rfloor} \begin{bmatrix} n \\ k \end{bmatrix} (-s)^k Luc^*_{n-2k}(x,s,q) = x^n \qquad (2.4)$$

and

$$\sum_{k=0}^{\lfloor \frac{n+1}{2} \rfloor} \left( \begin{bmatrix} n \\ k \end{bmatrix} - \begin{bmatrix} n \\ k-1 \end{bmatrix} \right) (-s)^k Fib_{n+1-2k}(x,s,q) = x^n \qquad (2.5)$$

hold (cf. [8], (3.1) and (3.2)).
Notice that $Luc_0(x,s,q) = 2$ and $Luc^*_0(x,s,q) = 1$, whereas $Luc_n(x,s,q) = Luc^*_n(x,s,q)$ for $n > 0$.

Fortunately there do exist $q-$analogues of the recurrences (1.19) and (1.21) which possess many of the looked for properties.

**Definition 2.1**

*The $q$-Chebyshev polynomials of the first kind are defined by the recurrence*

$$T_n(x,s,q) = (1+q^{n-1})xT_{n-1}(x,s,q) + q^{n-1}sT_{n-2}(x,s,q) \tag{2.6}$$

*with initial values $T_0(x,s,q) = 1$ and $T_1(x,s,q) = x$.*

The first terms are $1, x, [2]x^2 + qs, [4]x^3 + q[3]sx, \cdots$.

Some simple $q$-analogues of $T_n(1) = 1$ are

$$T_n(1,-1,q) = 1, \tag{2.7}$$

$$T_n\left(1,-\frac{1}{q},q\right) = q^{\binom{n}{2}}, \tag{2.8}$$

$$T_n(1,-q,q) = q^{\binom{n}{2}} + (1-q^n)\sum_{k=0}^{n-2} q^{\binom{k+1}{2}} \tag{2.9}$$

and

$$T_n(1,-q^2,q) = [n] - q^{n+1}[n-1]. \tag{2.10}$$

It is easily verified that

$$T_n\left(x,s,\frac{1}{q}\right) = \frac{T_n\left(x,\frac{s}{q},q\right)}{q^{\binom{n}{2}}}. \tag{2.11}$$

For $q = -1$ we get $T_{2n}(x,s,-1) = -sT_{2n-2}(x,s,-1)$ and
$T_{2n+1}(x,s,-1) = 2xT_{2n}(x,s,-1) + sT_{2n-1}(x,s,-1)$.
This gives the trivial sequence $(T_n(x,s,-q))_{n\geq 0} = (1, x, -s, -xs, s^2, s^2x, -s^3, -xs^3, \cdots)$. This is the reason for excluding $q = -1$.

**Proposition 2.1**

*The $q-$ Chebyshev polynomials of the first kind satisfy*

$$T_n(x,s,q) = \det\begin{pmatrix} x & qs & 0 & \cdots & 0 & 0 \\ -1 & (1+q)x & q^2s & \cdots & 0 & 0 \\ 0 & -1 & (1+q^2)x & \cdots & 0 & 0 \\ \vdots & \vdots & \vdots & \vdots & \vdots & \vdots \\ 0 & 0 & 0 & \cdots & (1+q^{n-2})x & q^{n-1}s \\ 0 & 0 & 0 & \cdots & -1 & (1+q^{n-1})x \end{pmatrix}.$$

This is easily seen by expanding this determinant with respect to the last column.

**Definition 2.2**

*The $q-$ Chebyshev polynomials of the second kind are defined by the recurrence*

$$U_n(x,s,q) = (1+q^n)xU_{n-1}(x,s,q) + q^{n-1}sU_{n-2}(x,s,q) \tag{2.12}$$

with initial values $U_0(x,s,q) = 1$ and $U_{-1}(x,s,q) = 0$.

The first terms are $1, [2]x, [4]x^2 + qs, [4](1+q^3)x^3 + q[4]sx, \cdots$.

Some simple $q-$ analogues of (1.4) are

$$U_n\left(1, -\frac{1}{q}, q\right) = q^{\binom{n}{2}}[n+1], \tag{2.13}$$

$$U_n(1,-1,q) = q^{\binom{n+1}{2}}\sum_{k=0}^{n}\frac{1}{q^{\binom{k+1}{2}}}. \tag{2.14}$$

$$U_n(1,-q,q) = \sum_{k=0}^{n} q^{\binom{k+1}{2}} \tag{2.15}$$

and

$$U_n(1,-q^2,q) = [n+1]. \tag{2.16}$$

It is easily verified that

$$U_n\left(s, x, \frac{1}{q}\right) = \frac{U_n(x,qs,q)}{q^{\binom{n+1}{2}}}. \tag{2.17}$$

For $q = -1$ we would have $U_{2n}(x,s,-1) = (-s)^n$ and $U_{2n+1}(x,s,-1) = 0$.

**Proposition 2.2**

*The $q$ – Chebyshev polynomials of the second kind satisfy*

$$U_n(x,s,q) = \det \begin{pmatrix} (1+q)x & qs & 0 & \cdots & 0 & 0 \\ -1 & (1+q^2)x & q^2 s & \cdots & 0 & 0 \\ 0 & -1 & (1+q^3)x & \cdots & 0 & 0 \\ \vdots & \vdots & \vdots & \vdots & \vdots & \vdots \\ 0 & 0 & 0 & \cdots & (1+q^{n-1})x & q^{n-1} s \\ 0 & 0 & 0 & \cdots & -1 & (1+q^n)x \end{pmatrix}.$$

In [3] and [16] a tiling interpretation of the classical Chebyshev polynomials has been given. This can easily be extended to the $q$ – case.

As in the classical case it is easier to begin with polynomials of the second kind.

We consider an $n \times 1$ – rectangle (called $n$ – board) where the $n$ cells of the board are numbered 1 to $n$. As in [3] and [16] we consider *tilings* with two sorts of squares, say white and black squares, and dominoes (which cover two adjacent cells of the board).

**Definition 2.3**

*To each tiling of a board we assign a weight $w$ in the following way: Each white square has weight $x$. A black square at position $i$ has weight $q^i x$ and a domino which covers positions $i-1, i$ has weight $q^{i-1} s$. The weight of a tiling is the product of its elements.*
*The weight of a set of tilings is the sum of their weights.*

Each tiling can be represented by a word in the letters $\{a, b, dd\}$. Here $a$ denotes a white square, $b$ a black square and $dd$ a domino.
For example the word $abbddaddaab$ represents the tiling with white squares at positions $1,6,9,10,$ black squares at $2,3,11$ and dominoes at $\{4,5\}$ and $\{7,8\}$. Its weight is
$x \cdot q^2 x \cdot q^3 x \cdot q^4 s \cdot x \cdot q^7 s \cdot x \cdot x \cdot q^{11} x = q^{27} s^2 x^7$.

**Theorem 2.1**

*The weight $w(V_n)$ of the set $V_n$ of all tilings of an $n$ – board is $w(V_n) = U_n(x,s,q)$.*

**Proof**
This holds for $n=1$ and $n=2$. Each $n$ – tiling $u_n$ has one of the following forms:
$u_{n-1}a, u_{n-1}b, u_{n-2}dd$.
Therefore

$$w(V_n) = \sum_{u_n \in V_n} w(u_n) = \sum_{u_{n-1} \in V_{n-1}} w(u_{n-1})x + \sum_{u_{n-1} \in V_{n-1}} w(u_{n-1})q^n x + \sum_{u_{n-2} \in V_{n-2}} w(u_{n-2})q^{n-1}s$$
$$= w(V_{n-1})(1+q^n)x + w(V_{n-2})q^{n-1}s$$

which implies Theorem 2.1.

**Remark 2.1**

If we more generally consider the weight $w_r$ which coincides with $w$ except that a black square at position $i$ has weight $q^i rx$ we get in the same way that $U_n^{(r)}(x,s,q) = w_r(V_n)$ satisfies

$$U_n^{(r)}(x,s,q) = (1+q^n r)xU_{n-1}^{(r)}(x,s,q) + q^{n-1}sU_{n-2}^{(r)}(x,s,q)$$

with initial values $U_0^{(r)}(x,s,q) = 1$ and $U_1^{(r)}(x,s,q) = (1+qr)x.$

In this case we get more generally

$$U_{m+n}^{(r)}(x,s,q) = U_m^{(r)}(x,s,q)U_n^{(q^m r)}(x,q^m s,q) + q^m sU_{m-1}^{(r)}(x,s,q)U_{n-1}^{(q^{m+1}r)}(x,q^{m+1}s,q).$$

The second term occurs when positions $(m, m+1)$ are covered by a domino and the first term in the other cases.

The same reasoning as above gives

**Proposition 2.3**

Let $u(n,k,s,r)$ be the $w_r$ – weight of all tilings on $\{1,\cdots,n\}$ with exactly $k$ dominoes. Then

$$u(n,k,s,r) = u(n-1,k,s,r)(1+q^n r)x + u(n-2,k-1,s,r)q^{n-1}s \qquad (2.18)$$

with initial values
$$u(n,0,s,r) = (1+qr)(1+q^2 r)\cdots(1+q^n r)x^n,$$

$u(1,0,s,r) = (1+qr)x$ and $u(1,k,s,r) = 0$ for $k>0.$

It is now easy to verify

**Theorem 2.2**

The $w_r$ – weight $u(n,k,s,r)$ of the set of all tilings on $\{1,\cdots,n\}$ with exactly $k$ dominoes is

$$u(n,k,s,r) = q^{k^2}\begin{bmatrix} n-k \\ k \end{bmatrix}(1+q^{k+1}r)\cdots(1+q^{n-k}r)s^k x^{n-2k} \qquad (2.19)$$

for $0 \le k \le \left\lfloor \dfrac{n}{2} \right\rfloor$ and $u(n,k,s,r) = 0$ for $k > \left\lfloor \dfrac{n}{2} \right\rfloor.$

**Proof**

The initial values coincide and by induction

$$u(n-1,k,s,r)(1+q^n r)x + u(n-2,k-1,s,r)q^{n-1}s = q^{k^2}\begin{bmatrix} n-k-1 \\ k \end{bmatrix}(1+q^{k+1}r)\cdots(1+q^{n-k-1}r)(1+q^n r)s^k x^{n-2k}$$

$$+q^{(k-1)^2}\begin{bmatrix} n-k-1 \\ k-1 \end{bmatrix}(1+q^k r)\cdots(1+q^{n-k-1}r)q^{n-1}s^k x^{n-2k}$$

$$= q^{k^2}(1+q^{k+1}r)\cdots(1+q^{n-k-1}r)s^k x^{n-2k}\left(\begin{bmatrix} n-k-1 \\ k \end{bmatrix}(1+q^n r) + \begin{bmatrix} n-k-1 \\ k-1 \end{bmatrix}q^{n-2k}(1+q^k r)\right)$$

$$= q^{k^2}(1+q^{k+1}r)\cdots(1+q^{n-k-1}r)s^k x^{n-2k}\left(\left(\begin{bmatrix} n-k-1 \\ k \end{bmatrix} + q^{n-2k}\begin{bmatrix} n-k-1 \\ k-1 \end{bmatrix}\right) + q^{n-k}r\left(q^k\begin{bmatrix} n-k-1 \\ k \end{bmatrix} + \begin{bmatrix} n-k-1 \\ k-1 \end{bmatrix}\right)\right)$$

$$= q^{k^2}(1+q^{k+1}r)\cdots(1+q^{n-k-1}r)s^k x^{n-2k}\left(\begin{bmatrix} n-k \\ k \end{bmatrix}(1+q^{n-k}r)\right).$$

Here we used the recurrence relations (2.1) for the $q$-binomial coefficients.

**Remark 2.2**

Formula (2.19) is the product of $q^{k^2}\begin{bmatrix} n-k \\ k \end{bmatrix}s^k x^{n-2k}$ and

$$(1+q^{k+1}r)\cdots(1+q^{n-k}r) = \sum_{\ell=0}^{n-2k}\begin{bmatrix} n-2k \\ \ell \end{bmatrix}(q^{(k+1)}r)^\ell q^{\binom{\ell}{2}}.$$

Ilse Fischer [12] has found a combinatorial reason for this product representation.

Let $v(n,k,\ell,x)$ be the $w_r$-weight of all tilings with $k$ dominoes and $\ell$ black squares. Then

$$v(n,k,\ell,x) = q^{k\ell}s^k r^\ell x^{n-2k}q^{k^2}\begin{bmatrix} n-k \\ k \end{bmatrix}q^{\binom{\ell+1}{2}}\begin{bmatrix} n-2k \\ \ell \end{bmatrix} = q^{k\ell}v(n,k,0,1)v(n-2k,0,\ell,x). \quad (2.20)$$

In order to give a combinatorial interpretation of this formula we observe that the weight can also be obtained from the following properties.

The fact that the weight of a domino at $\{i,i+1\}$ is $q^i s$ is equivalent with

a) each white square that appears before this domino contributes a $q$,
b) each black square that appears before this domino contributes a $q$,
c) each domino that appears before this domino contributes $q^2$
d) and the domino itself contributes $qs$.

The fact that the weight of a black square at $i$ is $q^i xr$ is equivalent with

e) each white square that appears before this black square contributes a $q$,
f) each black square that appears before this black square contributes a $q$,
g) each domino that appears before this black square contributes a $q^2$
h) and the black square itself contributes $qxr$.

This can also be reformulated in the following way:

1) Each black square contributes $qxr$,
2) each unordered pair of distinct black squares contributes a $q$,
3) each white square before a black square contributes a $q$,
4) each domino contributes $qs$,
5) each unordered pair of distinct dominoes contributes $q^2$,
6) each white square before a domino contributes a $q$,
7) each pair of a domino and a black square, where the order is irrelevant, contributes a $q$,
8) each domino before a black square contributes a $q$.

For b) and g) is equivalent with 7) and 8).

Now consider the right-hand side of (2.20).

Observe that $v(n,0,\ell,x)$ is determined by 1), 2) and 3); $v(n,k,0,x)$ is determined by 4), 5) and 6); and 7) gives $q^{k\ell}$.

We first distribute the dominoes on the $n-$board and let each unoccupied cell have weight 1. Then we distribute the white and black squares on the unoccupied cells. Their weight is $v(n-2k,0,\ell,x)$. The total weight of the configuration is $v(n,k,0,1)v(n-2k,0,\ell,x)$ if each black square before a domino contributes a $q$. For then 6) is satisfied for the computation of $v(n,k,0,1)$ since all squares contribute a $q$ (and thus behave as white squares in this context).

Thus the right-hand side of (2.20) satisfies 1) to 7), but instead of 8) we have
8´): each black square before a domino contributes a $q$.
Thus we must reverse the order of the dominoes and black squares to obtain (2.20).

An equivalent form is

**Proposition 2.4**

*Let $t$ be a tiling of an $n-$board with $k$ dominoes and $\ell$ black squares. Reverse the order of the dominoes and black squares in $t$ and obtain a tiling $T$. Denote by $A$ the tiling obtained by replacing in $T$ each square with a colourless square $c$ with weight $1$ and let $B$ be the tiling obtained by deleting all dominoes of $T$.*
*Then*
$$w_r(t) = q^{k\ell} w_r(A) w_r(B). \qquad (2.21)$$

**Example**
Consider the tiling $t = abbddaddaab$ with $(n,k,\ell) = (11,2,3)$ weight $q^{27}s^2r^3x^7$.
Then $T = abddddabaab$ and $A = ccddddccccc$ with $w_r(A) = q^8s^2$ and $B = ababaab$ with $w_r(B) = ababaab = q^{13}r^3x^7$.
This gives $w_r(t) = q^{2\cdot 3}\left(q^8s^2\right)\left(q^{13}x^7\right)$.
Theorem 2.2 implies for $r = 1$

**Theorem 2.3**

$$U_n(x,s,q) = \sum_{k=0}^{\left\lfloor\frac{n}{2}\right\rfloor} q^{k^2} \begin{bmatrix} n-k \\ k \end{bmatrix} (1+q^{k+1})\cdots(1+q^{n-k})s^k x^{n-2k}. \quad (2.22)$$

For the $q$-Chebyshev polynomials of the first kind the situation is somewhat more complicated. Here we get

**Theorem 2.4**

$T_n(x,s,q)$ is the weight of the subset of all tilings of $\{1,\cdots,n\}$ where the last block is either a white square or a domino.
Therefore for $n > 0$

$$T_n(x,s,q) = xU_{n-1}(x,s,q) + q^{n-1}sU_{n-2}(x,s,q). \quad (2.23)$$

**Proof**
It suffices to prove that the right-hand side satisfies the initial values and the recurrence (2.6).

$(1+q^n)x\left(xU_{n-1}(x,s,q) + q^{n-1}sU_{n-2}(x,s,q)\right) + q^n s\left(xU_{n-2}(x,s,q) + q^{n-2}sU_{n-3}(x,s,q)\right)$
$= x^2 U_{n-1}(x,s,q) + q^n x^2 U_{n-1}(x,s,q) + q^{n-1}sxU_{n-2}(x,s,q) + q^{2n-1}sxU_{n-2}(x,s,q)$
$+ q^n sxU_{n-2}(x,s,q) + q^{2n-2}s^2 U_{n-3}(x,s,q)$
$= x\left((1+q^n)xU_{n-1}(x,s,q) + q^{n-1}sU_{n-2}(x,s,q)\right)$
$+ q^n s\left((1+q^{n-1})xU_{n-2}(x,s,q) + q^{n-2}sU_{n-4}(x,s,q)\right)$
$= xU_n(x,s,q) + q^n sU_{n-1}(x,s,q).$

**Theorem 2.5**

*The $q$ – Chebyshev polynomials of the first kind are given by*

$$T_n(x,s,q) = \sum_{k=0}^{\lfloor \frac{n}{2} \rfloor} q^{k^2} \frac{(1+q)\cdots(1+q^{n-1})}{(1+q)\cdots(1+q^k)\cdot(1+q^{n-k})\cdots(1+q^{n-1})} \frac{[n]}{[n-k]} \begin{bmatrix} n-k \\ k \end{bmatrix} s^k x^{n-2k}$$

(2.24)

$$= \sum_{k=0}^{\lfloor \frac{n-1}{2} \rfloor} q^{k^2} (1+q^{k+1})\cdots(1+q^{n-k-1}) \frac{[n]}{[n-k]} \begin{bmatrix} n-k \\ k \end{bmatrix} s^k x^{n-2k} + [n \equiv 0 \bmod 2] q^{n^2} s^n.$$

**Proof**

Consider the subset of all tilings of an $n$ – board whose last block is not a black square. Let $t(n,k,s)$ be the weight of all these tilings with exactly $k$ dominoes.
Then

$$t(n,k,s) = u(n-1,k,s)x + u(n-2,k-1)q^{n-1}s. \tag{2.25}$$

We show that

$$t(n,k,s) = q^{k^2} \frac{(1+q)\cdots(1+q^{n-1})}{(1+q)\cdots(1+q^k)\cdot(1+q^{n-k})\cdots(1+q^{n-1})} \frac{[n]}{[n-k]} \begin{bmatrix} n-k \\ k \end{bmatrix} s^k x^{n-2k}. \tag{2.26}$$

This is true for $n=1$ and $n=2$. By induction we get for $2k \leq n-1$

$$t(n,k,s) = u(n-1,k,s)x + u(n-2,k-1)q^{n-1}s$$

$$= q^{k^2} \begin{bmatrix} n-k-1 \\ k \end{bmatrix} (1+q^{k+1})\cdots(1+q^{n-k-1}) s^k x^{n-2k}$$

$$+ q^{n-1} q^{k^2-2k+1} \begin{bmatrix} n-1-k \\ k-1 \end{bmatrix} (1+q^k)\cdots(1+q^{n-1-k}) s^k x^{n-2k}$$

$$= q^{k^2}(1+q^{k+1})\cdots(1+q^{n-k-1}) \left( \begin{bmatrix} n-k-1 \\ k \end{bmatrix} + (1+q^k)q^{n-2k} \begin{bmatrix} n-1-k \\ k-1 \end{bmatrix} \right) s^k x^{n-2k}$$

$$= q^{k^2}(1+q^{k+1})\cdots(1+q^{n-k-1}) \left( \begin{bmatrix} n-k \\ k \end{bmatrix} + q^{n-k} \begin{bmatrix} n-1-k \\ k-1 \end{bmatrix} \right) s^k x^{n-2k}$$

$$= q^{k^2}(1+q^{k+1})\cdots(1+q^{n-k-1}) \frac{[n]}{[n-k]} \begin{bmatrix} n-k \\ k \end{bmatrix} s^k x^{n-2k}$$

and for $2k = 2n$

$$t(2n,n,s) = u(2n-1,n,s)x + u(2n-2,n-1)q^{2n-1}s$$

$$= q^{2n-1} q^{n^2-2n+1} \begin{bmatrix} n-1 \\ n-1 \end{bmatrix} s^n = q^{n^2} s^n.$$

For $q=1$ the polynomial $T_n(x,s)$ can also be interpreted as the weight of the set $T_n$ of all tilings which begin with a domino or with a white square since in this case the weights of the words $c_1 \cdots c_n$ and $c_n \cdots c_1$ coincide.

In the general case this is not true. For example for $n=2$ the set $T_2 = \{aa, ab, dd\}$ has weight $w(T_2) = x^2 + q^2 x^2 + qs \neq T_2(x,s,q) = x^2 + qx^2 + qs$.

But we have

**Theorem 2.6**

$$T_n(x,s,q) = xU_{n-1}(x, q^2 s, q) + qsU_{n-2}(x, q^2 s, q). \tag{2.27}$$

**Proof**

It suffices to show that the right-hand side satisfies recurrence (2.6).

$$\left(1+q^{n-1}\right)x\left(xU_{n-2}(x,q^2s,q) + qsU_{n-3}(x,q^2s,q)\right) + q^{n-1}s\left(xU_{n-3}(x,q^2s,q) + qsU_{n-4}(x,q^2s,q)\right)$$
$$= x^2 U_{n-2}(x,q^2s,q) + q^{n-1} x^2 U_{n-2}(x,q^2s,q) + qsxU_{n-3}(x,q^2s,q) + q^n sxU_{n-3}(x,q^2s,q)$$
$$+ q^{n-1} sxU_{n-3}(x,q^2s,q) + q^n s^2 U_{n-4}(x,q^2s,q)$$
$$= x\left(\left(1+q^{n-1}\right)xU_{n-2}(x,q^2s,q) + q^{n-2} q^2 s U_{n-3}(x,q^2s,q)\right)$$
$$+ qs\left(\left(1+q^{n-2}\right)xU_{n-3}(x,q^2s,q) + q^{n-3} q^2 s U_{n-4}(x,q^2s,q)\right)$$
$$= xU_{n-1}(x,q^2s,q) + qsU_{n-2}(x,q^2s,q).$$

Theorem 2.6 has the following tiling interpretation:

Define another weight $W$ such that each white square has weight $x$, each black square at position $i$ has weight $q^i x$ and each domino at position $(i-1,i)$ has weight $q^{i+1}s$ if $i<n$. But a domino at position $(n-1,n)$ has weight $qs$.
If we join the ends of the board to a circle such that the position after $n$ is $1$ this can also be formulated as: If $(i-1,i,j)$ are consecutive points then a domino at position $(i-1,i)$ has weight $q^j s$. Then $T_n(x,s,q)$ is the weight of all such tilings which have no black square at position $n$. (Note that on the circle there are no dominoes at position $(n,1)$.)

In order to find a $q-$analogue of (1.23) let us first consider this identity in more detail.
$$\left(x + \sqrt{x^2+s}\right)^n = T_n(x,s) + U_{n-1}(x,s)\sqrt{x^2+s}$$
is equivalent with
$$T_{n+1}(x,s) + U_n(x,s)\sqrt{x^2+s} = \left(x+\sqrt{x^2+s}\right)^{n+1} = \left(x+\sqrt{x^2+s}\right)\left(x+\sqrt{x^2+s}\right)^n$$
$$= \left(x+\sqrt{x^2+s}\right)\left(T_n(x,s) + U_{n-1}(x,s)\sqrt{x^2+s}\right)$$
$$= T_n(x,s)x + \left(x^2+s\right)U_{n-1}(x,s) + \left(T_n(x,s) + U_{n-1}(x,s)x\right)\sqrt{x^2+s}.$$

Therefore (1.23) is equivalent with both identities

$$T_{n+1}(x,s) = T_n(x,s)x + (x^2 + s)U_{n-1}(x,s) \tag{2.28}$$

and

$$U_n(x,s) = T_n(x,s) + U_{n-1}(x,s)x. \tag{2.29}$$

To prove identity (2.28) observe that for $q = 1$ a tiling of an $(n+1)$ – board which does not end with a black square either ends with two white squares *aa* or with a domino and a white square *dda*. The weight $w$ of these tilings is $T_n(x,s)x$. Or it ends with *ba* or *dd*. Their weight is $(x^2 + s)U_{n-1}(x,s)$.

Identity (2.29) simply means that an arbitrary tiling either ends with a black square which gives the weight $U_{n-1}(x,s)x$ or does not end with a black square which gives $T_n(x,s)$.

For arbitrary $q$ this classification of the tilings implies the identities

$$T_{n+1}(x,s,q) = xT_n(x,s,q) + q^n(x^2 + s)U_{n-1}(x,s,q) \tag{2.30}$$

and

$$U_n(x,s,q) = T_n(x,s,q) + q^n x U_{n-1}(x,s,q). \tag{2.31}$$

But there is also another $q$ – analogue of (2.28):

$$T_{n+1}(x,s,q) = q^n x T_n(x,s,q) + (x^2 + qs)U_{n-1}(x,q^2 s,q). \tag{2.32}$$

By (2.27) we have $T_{n+1}(x,s,q) = xU_n(x,q^2 s,q) + qsU_{n-1}(x,q^2 s,q)$.
Therefore by (2.12)

$$U_n(x,q^2 s,q) - xU_{n-1}(x,q^2 s,q) = q^n x U_{n-1}(x,q^2 s,q) + q^{n+1} s U_{n-2}(x,q^2 s,q)$$
$$= q^n \left( xU_{n-1}(x,q^2 s,q) + qs U_{n-2}(x,q^2 s,q) \right) = q^n T_n(x,s,q).$$

Thus

$$U_n(x,q^2 s,q) = q^n T_n(x,s,q) + xU_{n-1}(x,q^2 s,q) \tag{2.33}$$

and (2.27) implies (2.32).

As $q$ – analogue of (2.28) and (2.29) we can now choose the identities (2.31) and (2.32) which we write in the form

$$\begin{aligned} T_{n+1}(x,s,q) &= q^n x T_n(x,s,q) + (x^2 + qs)\eta^2 U_{n-1}(x,s,q) \\ U_n(x,s,q) &= T_n(x,s,q) + q^n x U_{n-1}(x,s,q). \end{aligned} \tag{2.34}$$

Here $\eta$ denotes the linear operator on the polynomials in $s$ defined by $\eta p(s) = p(qs)$.

To stress the analogy with (1.23) we introduce a formal square root $A = \sqrt{(x^2+s)\eta^2}$ which commutes with $x$ and real or complex numbers and satisfies $A^2 = (x^2+qs)\eta^2$ and write (2.34) in the form

$$T_{n+1}(x,s,q) + AU_n(x,s,q) = (q^n x + A)(T_n(x,s,q) + AU_{n-1}(x,s,q)). \tag{2.35}$$

Since $(q^i x + A)(q^j x + A) = (q^j x + A)(q^i x + A)$ using the $q$-binomial theorem (2.3) we get as analogue of (1.23)

$$p_n(x,A) = (x+A)(qx+A)\cdots(q^{n-1}x+A) = T_n(x,s,q) + AU_{n-1}(x,s,q). \tag{2.36}$$

This gives

**Theorem 2.7**

*For the $q$-Chebyshev polynomials the following formulae hold:*

$$T_n(x,s,q) = \frac{p_n(x,A) + p_n(x,-A)}{2}\mathbf{1} = \sum_{k=0}^{\lfloor n/2 \rfloor} q^{\binom{n-2k}{2}} \begin{bmatrix} n \\ 2k \end{bmatrix} x^{n-2k} \prod_{j=0}^{k-1}(x^2 + q^{2j+1}s) \tag{2.37}$$

*and*

$$U_n(x,s,q) = \frac{p_{n+1}(x,A) - p_{n+1}(x,-A)}{2A}\mathbf{1} = \sum_{k=0}^{\lfloor n/2 \rfloor} q^{\binom{n-2k}{2}} \begin{bmatrix} n+1 \\ 2k+1 \end{bmatrix} x^{n-2k} \prod_{j=0}^{k-1}(x^2 + q^{2j+1}s). \tag{2.38}$$

**Proof**

This follows from (2.3) and the observation that $A^{2k} = \left((x^2+qs)\eta^2\right)^k = \prod_{j=0}^{k-1}(x^2+q^{2j+1}s)\eta^{2k}$.

If we expand $\prod_{j=0}^{k-1}(x^2 + q^{2j+1}s) = \sum_{j=0}^{k} q^{j^2} s^j \begin{bmatrix} k \\ j \end{bmatrix}_{q^2} x^{2k-2j}$

we get by comparing coefficients in (2.37) and (2.38)

**Theorem 2.8**

*For $j \leq n$ the identities*

$$\sum_{k=0}^{\lfloor n/2 \rfloor} q^{\binom{n-2k}{2}} \begin{bmatrix} n \\ 2k \end{bmatrix}_q \begin{bmatrix} k \\ j \end{bmatrix}_{q^2} = \frac{(1+q)\cdots(1+q^{n-1})}{(1+q)\cdots(1+q^j)\cdot(1+q^{n-j})\cdots(1+q^{n-1})} \frac{[n]_q}{[n-j]_q} \begin{bmatrix} n-j \\ j \end{bmatrix}_q \tag{2.39}$$

*and*

$$\sum_{k=0}^{\left\lfloor\frac{n}{2}\right\rfloor} q^{\binom{n-2k}{2}} \begin{bmatrix} n+1 \\ 2k+1 \end{bmatrix}_q \begin{bmatrix} k \\ j \end{bmatrix}_{q^2} = (1+q^{j+1})\cdots(1+q^{n-j}) \begin{bmatrix} n-j \\ j \end{bmatrix}_q \tag{2.40}$$

hold.

## Remark 2.3

It would be nice to find a combinatorial interpretation of these identities.
For $q=1$ we get from (1.23)
$T_n(x,s)^2 - (x^2+s)U_{n-1}(x,s)^2 = (-s)^n.$
Since $A$ does not commute with polynomials in $s$ we cannot deduce a $q-$analogue of this formula from (2.36).

But we can instead consider the matrices

$$A_n = \begin{pmatrix} x & q^n(x^2+s) \\ 1 & q^n x \end{pmatrix}. \tag{2.41}$$

We then get

**Theorem 2.9**

$$\begin{pmatrix} T_n(x,s,q) & (x^2+s)U_{n-1}(x,qs,q) \\ U_{n-1}(x,s,q) & T_n\left(x,\dfrac{s}{q},q\right) \end{pmatrix} = A_{n-1}A_{n-2}\cdots A_0. \tag{2.42}$$

**Proof**

We must show that

$$\begin{pmatrix} T_{n+1}(x,s,q) & (x^2+s)U_n(x,qs,q) \\ U_n(x,s,q) & T_{n+1}\left(x,\dfrac{s}{q},q\right) \end{pmatrix} = \begin{pmatrix} x & q^n(x^2+s) \\ 1 & q^n x \end{pmatrix} \begin{pmatrix} T_n(x,s,q) & (x^2+s)U_{n-1}(x,qs,q) \\ U_{n-1}(x,s,q) & T_n\left(x,\dfrac{s}{q},q\right) \end{pmatrix}$$

or equivalently

$T_{n+1}(x,s,q) = xT_n(x,s,q) + q^n(x^2+s)U_{n-1}(x,s,q),$

$U_n(x,s,q) = T_n(x,s,q) + q^n x U_{n-1}(x,s,q),$

$U_n(x,q^2 s,q) = q^n T_n(x,s,q) + xU_{n-1}(x,q^2 s,q),$

$T_{n+1}(x,s,q) = q^n x T_n(x,s,q) + (x^2+qs)U_{n-1}(x,q^2 s,q).$

This follows from the recurrences (2.30), (2.31), (2.32) and (2.33).

If we take determinants in (2.42) we get the desired $q-$analogue of
$T_n(x,s)^2 - (x^2+s)U_{n-1}(x,s)^2 = (-s)^n.$

**Theorem 2.10**

$$T_n(x,s,q)T_n(x,qs,q) - (x^2+qs)U_{n-1}(x,qs,q)U_{n-1}(x,q^2s,q) = q^{\binom{n+1}{2}}(-s)^n. \qquad (2.43)$$

For example for $(x,s) = (1,-1)$ this reduces to

$$T_n(1,-q,q) - (1-q)\sum_{k=1}^{n} q^{\binom{k}{2}}[n] = T_n(1,-q,q) - (1-q^n)\sum_{k=1}^{n} q^{\binom{k}{2}} = q^{\binom{n+1}{2}}.$$

In [11] many other identities occur. These follow in an easy manner from the identities obtained above.

Since the $q-$Chebyshev polynomials satisfy a three-term recurrence they are orthogonal with respect to some linear functionals, i.e. $L(T_n(x,s,q)T_m(x,s,q)) = 0$ and
$M(U_n(x,s,q)U_m(x,s,q)) = 0$ for $n \neq m$.
These linear functionals are uniquely determined by
$L(T_n(x,s,q)) = [n=0]$ and $M(U_n(x,s,q)) = [n=0]$.

These linear functionals are closely related. From (2.30) we get
$T_{n+1}(x,s,q) - xT_n(x,s,q) = q^n(x^2+s)U_{n-1}(x,s,q).$

By (2.6) we have $xT_n(x,s,q) = \dfrac{T_{n+1}(x,s,q) - q^n sT_{n-1}(x,s,q)}{1+q^n}$

and therefore we obtain

$$T_{n+1}(x,s,q) + sT_{n-1}(x,s,q) = (1+q^n)(x^2+s)U_{n-1}(x,s,q). \qquad (2.44)$$

If we apply the linear functional $L$ to this identity we deduce that

$$(1+q)L\left(\left(1+\frac{x^2}{s}\right)U_n(x,s,q)\right) = [n=0] = M(U_n(x,s,q)). \qquad (2.45)$$

By linearity we obtain

$$(1+q)L\left(\left(1+\frac{x^2}{s}\right)p(x)\right) = M(p(x)) \qquad (2.46)$$

for all polynomials $p(x)$.

As $q-$analogue of (1.14) we get

$$L\left(T_n^2\right) = \begin{cases} 1 & \text{if } n=0 \\ \dfrac{q^{\binom{n+1}{2}}(-s)^n}{1+q^n} & \text{if } n>0 \end{cases} \tag{2.47}$$

This follows by applying $L$ to (2.6) which gives $L\left(x^n T_n\right) = -\dfrac{q^n s}{1+q^n} L\left(x^{n-1} T_{n-1}\right)$ and therefore

$$L\left(x^n T_n\right) = (-s)^n \frac{q^{\binom{n+1}{2}}}{(1+q)(1+q^2)\cdots(1+q^n)}.$$

Now observe that $L\left(T_n^2\right) = L\left((1+q)\cdots(1+q^{n-1}) x^n T_n\right)$.

Of special interest are the moments of these linear functionals, i.e. the values $L(x^n)$ and $M(x^n)$. To find these values it suffices to find the uniquely determined representation of $x^n$ as a linear combination of the $q-$Chebyshev polynomials.
These have been calculated in [11] for the corresponding monic polynomials. Therefore I only state the results in the present notation:

For the $q-$Chebyshev polynomials of the first kind we have

$$x^n = \sum_{k=0}^{\lfloor n/2 \rfloor} \begin{bmatrix} n \\ k \end{bmatrix} (1+q^{n-2k}[2k\neq n])(-qs)^k \frac{T_{n-2k}(x,s,q)}{(1+q)\cdots(1+q^k)(1+q)\cdots(1+q^{n-k})}. \tag{2.48}$$

This gives as $q-$analogue of (1.15)

$$L\left(x^{2n}\right) = \begin{bmatrix} 2n \\ n \end{bmatrix} \frac{(-qs)^n}{\prod_{j=1}^{n}(1+q^j)^2} \tag{2.49}$$

and $L(x^{2n+1}) = 0$.

For the monic polynomials we get the three-term recurrence with $s(n)=0$, $t(0) = \dfrac{qs}{1+q}$ and

$$t(n) = \frac{q^{n+1} s}{(1+q^n)(1+q^{n+1})}.$$

For the $q$-Chebyshev polynomials of the second kind the corresponding formulae are

$$M(U_n^2) = (-s)^n q^{\binom{n+1}{2}} \frac{1+q}{1+q^{n+1}} \qquad (2.50)$$

as $q$-analogue of (1.17) and

$$x^n = \sum_{k=0}^{\left\lfloor \frac{n}{2} \right\rfloor} \left( \begin{bmatrix} n \\ k \end{bmatrix} - \begin{bmatrix} n \\ k-1 \end{bmatrix} \right)(-s)^k \frac{1+q^{n-2k+1}}{\prod_{j=1}^{k}(1+q^j) \prod_{j=1}^{n-k+1}(1+q^j)} U_{n-2k}(x,s,q) \qquad (2.51)$$

and therefore

$$M(x^{2n}) = \frac{1}{[n+1]} \begin{bmatrix} 2n \\ n \end{bmatrix} \frac{1+q}{1+q^{n+1}} \frac{(-qs)^n}{\prod_{j=1}^{n}(1+q^j)^2} \qquad (2.52)$$

and $M(x^{2n+1}) = 0$.
Of course (2.52) also follows directly from (2.49) and (2.46).

The parameters for the three-term recurrence of the monic polynomials are $s(n) = 0$ and
$$t(n) = \frac{q^{n+1}s}{(1+q^{n+1})(1+q^{n+2})}.$$

**Remark 2.4**

The $q$-Chebyshev polynomials have also appeared, partly implicitly and without recognizing them as $q$-analogues of the Chebyshev polynomials, in [6], [7] and [13] in the course of computing Hankel determinants of $\mu_n = \frac{(aq;q)_n}{(abq^2;q)_n}$, which are the moments of the little $q$-Jacobi polynomials $p_n(x;a,b\,|\,q)$ (cf. [14]). Note that $L(x^{2n}) = \frac{(q;q^2)_n}{(q^2;q^2)_n}(-qs)^n$ and

$$M(x^{2n}) = \frac{(q^2;q^2)_n}{(q^4;q^2)_n}(-qs)^n.$$

## 3. Some further properties

The $q$-Chebyshev polynomials $T_{2n}(1,s,q), T_{2n+1}(1,s,q), U_{2n}(1,s,q)$ and $U_{2n+1}(1,s,q)$ are polynomials in $s$ of degree $n$.
Therefore there exist unique representations

$$T_{2n+1}(1,s,q) = \sum_{k=0}^{n} a(n,k,q) T_{2k}(1,s,q) \qquad (3.1)$$

and

$$U_{2n+1}(1,s,q) = \sum_{k=0}^{n} b(n,k,q) U_{2k}(1,s,q). \tag{3.2}$$

To obtain these representations we need $q-$analogues of the tangent and Genocchi numbers. The $q-tangent\ numbers$ $t_{2n+1}(q)$ are well-known objects defined by the generating function

$$\frac{e(z)-e(-z)}{e(z)+e(-z)} = \sum_{n\geq 0} \frac{(-1)^n t_{2n+1}(q)}{[2n+1]!} z^{2n+1}. \tag{3.3}$$

**Theorem 3.1**

$$T_{2n+1}(x,s,q) = \sum_{k=0}^{n} \begin{bmatrix} 2n+1 \\ 2k \end{bmatrix} (-1)^{n-k} t_{2n-2k+1}(q) x^{2n+1-2k} T_{2k}(x,s,q). \tag{3.4}$$

**Proof**

In (2.37) we have seen that $T_n(1,s,q) = \sum_{k=0}^{\lfloor \frac{n}{2} \rfloor} \begin{bmatrix} n \\ 2k \end{bmatrix} q^{\binom{n-2k}{2}} (1+qs)(1+q^3s)\cdots(1+q^{2k-1}s).$

This implies that

$$T(z,s,q) = \sum_{n\geq 0} \frac{T_n(1,s,q)}{[n]!} z^n \tag{3.5}$$

satisfies

$$T(z,s,q) = \frac{1}{e(-z)} \sum_{n\geq 0} (1+qs)(1+q^3s)\cdots(1+q^{2n-1}s) \frac{z^{2n}}{[2n]!}. \tag{3.6}$$

Therefore $e(-z)T(z,s,q) = e(z)T(-z,s,q)$ and
$$\big(e(z)-e(-z)\big)\big(T(z,s,q)+T(-z,s,q)\big) = \big(e(z)+e(-z)\big)\big(T(z,s,q)-T(-z,s,q)\big)$$
or

$$\frac{\sum_{n\geq 0} \frac{T_{2n+1}(1,s,q)}{[2n+1]!} z^{2n+1}}{\sum_{n\geq 0} \frac{T_{2n}(1,s,q)}{[2n]!} z^{2n}} = \frac{e(z)-e(-z)}{e(z)+e(-z)} = \sum_{n\geq 0} \frac{(-1)^n t_{2n+1}(q)}{[2n+1]!} z^{2n+1}. \tag{3.7}$$

Note that the left-hand side does not depend on $s$. If we choose $s=0$ we get that

$$\frac{\sum_{n\geq 0} \frac{(-q;q)_{2n}}{[2n+1]!} z^{2n+1}}{1+\sum_{n\geq 1} \frac{(-q;q)_{2n-1}}{[2n]!} z^{2n}} = \frac{e(z)-e(-z)}{e(z)+e(-z)}. \tag{3.8}$$

(3.7) implies

$$\sum_{n\geq 0}\frac{T_{2n+1}(1,s,q)}{[2n+1]!}z^{2n+1} = \sum_{n\geq 0}\frac{(-1)^n t_{2n+1}(q)}{[2n+1]!}z^{2n+1}\sum_{n\geq 0}\frac{T_{2n}(1,s,q)}{[2n]!}z^{2n}$$

which gives by comparing coefficients

$$T_{2n+1}(1,s,q) = \sum_{k=0}^{n}\begin{bmatrix}2n+1\\2k\end{bmatrix}(-1)^{n-k}t_{2n-2k+1}(q)T_{2k}(1,s,q) \tag{3.9}$$

and therefore also (3.4).

For $q=1$ the Chebyshev polynomials satisfy

$$\sum_{j=0}^{n}\binom{n}{j}(-2x)^j T_{2n+m-j}(x,s) = s^n T_m(x,s) \tag{3.10}$$

and

$$\sum_{j=0}^{n}\binom{n}{j}(-2x)^j U_{2n+m-1-j}(x,s) = s^n U_m(x,s). \tag{3.11}$$

For these identities are equivalent with

$$\sum_{j=0}^{n}\binom{n}{j}(-2x)^j \left(x+\sqrt{x^2+s}\right)^{2n+m-j} = s^n \left(x+\sqrt{x^2+s}\right)^{2n+m-j}$$

which in turn reduces to the trivial identity

$$\left(x+\sqrt{x^2+s}\right)^{n+m}\left(x+\sqrt{x^2+s}-2x\right)^n = \left(x+\sqrt{x^2+s}\right)^m\left(\sqrt{x^2+s}+x\right)^n\left(\sqrt{x^2+s}-x\right)^n = s^n\left(x+\sqrt{x^2+s}\right)^m.$$

In order to simplify the exposition we let $x=1$ and prove as $q-$analogue of (3.10)

**Theorem 3.2**

$$\sum_{j=0}^{n}(-1)^j q^{\binom{j}{2}}\begin{bmatrix}n\\j\end{bmatrix}\prod_{i=n+m+1-j}^{n+m}(1+q^i)T_{2n+m-j}(1,s,q) = q^{n^2+mn}s^n T_m(1,s,q). \tag{3.12}$$

**Proof**

Let $m\in\mathbb{N}$. We consider the following matrix $\left(a(n,k,m)\right)_{n,k\geq 0}$ with $a(n,k,m) = s^k T_{n-k+m}(1,s,q)$ for $0\leq k\leq n$ and $a(n,k,m) = 0$ for $k>n$. The first terms are

$$\begin{pmatrix} T_m(1,s,q) & & & & \\ T_{m+1}(1,s,q) & sT_m(1,s,q) & & & \\ T_{m+2}(1,s,q) & sT_{m+1}(1,s,q) & s^2T_m(1,s,q) & & \\ T_{m+3}(1,s,q) & sT_{m+2}(1,s,q) & s^2T_{m+1}(1,s,q) & s^3T_m(1,s,q) & \\ \vdots & \vdots & \vdots & \vdots & \ddots \end{pmatrix}.$$

The recurrence for $T_n(1,s,q)$ gives

$$a(n,k,m) = \frac{a(n+1,k-1,m)-(1+q^{n+m+1-k})a(n,k-1,m)}{q^{n+m+1-k}}.$$

This implies that

$$a(n,k,m) = \frac{1}{q^{k(n+m)}} \sum_{j=0}^{k} (-1)^j q^{\binom{j}{2}} \begin{bmatrix} k \\ j \end{bmatrix} \prod_{i=n+m+1-j}^{n+m} (1+q^i) T_{n+m+k-j}(1,s,q).$$

This is true for $k=0$.
If it holds for $k-1$ then

$$a(n,k,m) = \frac{a(n+1,k-1,m)-(1+q^{n+m-1-k})a(n,k-1,m)}{q^{n+m+1-k}}$$

$$= \frac{1}{q^{n+m+1-k}} \frac{1}{q^{(k-1)(n+m+1)}} \sum_{j=0}^{k-1} (-1)^j q^{\binom{j}{2}} \begin{bmatrix} k-1 \\ j \end{bmatrix} \prod_{i=n+m+2-j}^{n+m+1} (1+q^i) T_{n+m+k-j}(1,s,q)$$

$$-(1+q^{n+m+1-k}) \frac{1}{q^{n+m+1-k}} \frac{1}{q^{(k-1)(n+m)}} \sum_{j=0}^{k-1} (-1)^j q^{\binom{j}{2}} \begin{bmatrix} k-1 \\ j \end{bmatrix} \prod_{i=n+m+1-j}^{n+m} (1+q^i) T_{n+m+k-1-j}(1,s,q)$$

$$= \frac{1}{q^{k(n+m)}} \sum_{j=0}^{k-1} (-1)^j q^{\binom{j}{2}} \begin{bmatrix} k-1 \\ j \end{bmatrix} \prod_{i=n+m+2-j}^{n+m+1} (1+q^i) T_{n+m+k-j}(1,s,q)$$

$$+(1+q^{n+m+1-k}) \frac{q^{k-1}}{q^{k(n+m)}} \sum_{j=1}^{k} (-1)^j q^{\binom{j-1}{2}} \begin{bmatrix} k-1 \\ j-1 \end{bmatrix} \prod_{i=n+m+2-j}^{n+m} (1+q^i) T_{n+m+k-j}(1,s,q)$$

$$= \frac{1}{q^{k(n+m)}} \sum_{j=0}^{k} (-1)^j q^{\binom{j}{2}} \prod_{i=n+m+2-j}^{n+m} (1+q^i) T_{n+m+k-j}(1,s,q) \left( \begin{bmatrix} k-1 \\ j \end{bmatrix}(1+q^{n+m+1}) + q^{k-j}(1+q^{n+m+1-k}) \begin{bmatrix} k-1 \\ j-1 \end{bmatrix} \right)$$

$$= \frac{1}{q^{k(n+m)}} \sum_{j=0}^{k} (-1)^j q^{\binom{j}{2}} \begin{bmatrix} k \\ j \end{bmatrix} \prod_{i=n+m+1-j}^{n+m} (1+q^i) T_{n+m+k-j}(1,s,q).$$

This gives (3.12).

As special cases we get for $m=0$ and $m=1$

$$\sum_{j=0}^{n} (-1)^j q^{\binom{j}{2}} \begin{bmatrix} n \\ j \end{bmatrix} \prod_{i=n+1-j}^{n} (1+q^i) T_{2n-j}(1,s,q) = q^{n^2} s^n$$

and

$$\sum_{j=0}^{n}(-1)^{j}q^{\binom{j}{2}}\begin{bmatrix}n\\j\end{bmatrix}\prod_{i=n+2-j}^{n+1}(1+q^{i})T_{2n+1-j}(1,s,q)=q^{n^{2}+n}s^{n}.$$

This implies

$$q^{n}\sum_{j=1}^{n+1}(-1)^{j-1}q^{\binom{j-1}{2}}\begin{bmatrix}n\\j-1\end{bmatrix}\prod_{i=n+2-j}^{n}(1+q^{i})T_{2n+1-j}(1,s,q)$$
$$=\sum_{j=1}^{n+1}(-1)^{j}q^{\binom{j}{2}}\begin{bmatrix}n\\j\end{bmatrix}\prod_{i=n+2-j}^{n+1}(1+q^{i})T_{2n+1-j}(1,s,q)+T_{2n+1}(1,s,q)$$

or

$$\sum_{j=1}^{n+1}(-1)^{j-1}q^{\binom{j}{2}}\prod_{i=n+2-j}^{n}(1+q^{i})\left(\begin{bmatrix}n+1\\j\end{bmatrix}+q^{n+1}\begin{bmatrix}n\\j\end{bmatrix}\right)T_{2n+1-j}(1,s,q)=T_{2n+1}(1,s,q). \quad (3.13)$$

Of course we could also replace $\begin{bmatrix}n+1\\j\end{bmatrix}+q^{n+1}\begin{bmatrix}n\\j\end{bmatrix}$ by $\begin{bmatrix}n+1\\2j\end{bmatrix}\frac{[2n+2-2j]}{[n+1]}$.

Define now a linear functional $\mu$ on the polynomials in $s$ by $\mu(T_{2n}(1,s,q))=[n=1]$. Then by (3.9) $\mu(T_{2n+1}(1,s,q))=(-1)^{n}t_{2n+1}(q)$.

Thus we get the following identities for the $q-$tangent numbers

$$t_{2n+1}(q)=\sum_{j=1}^{\lfloor\frac{n+1}{2}\rfloor}(-1)^{j-1}q^{\binom{2j}{2}}\prod_{i=n+2-2j}^{n}(1+q^{i})\begin{bmatrix}n+1\\2j\end{bmatrix}\frac{[2n+2-2j]}{[n+1]}t_{2n+1-2j}(q). \quad (3.14)$$

For $q=1$ this reduces to

$$t_{2n+1}=\sum_{j=1}^{\lfloor\frac{n+1}{2}\rfloor}(-1)^{j-1}2^{2j}\binom{n+1}{2j}\frac{n+1-j}{n+1}. \quad (3.15)$$

The first identities are

$t_{3}=2t_{1},\ t_{5}=8t_{3},\ t_{7}=18t_{5}-8t_{3},\ t_{9}=32t_{7}-48t_{5},\ t_{11}=50t_{9}-160t_{7}+32t_{5}.$

What at first glance appears as a new identity turns out to be an old acquaintance if we use (1.35) and write (3.15) in terms of Genocchi numbers. For then we get

$$\sum_{j=0}^{n}(-1)^{j}\binom{n}{2j}G_{2n-2j}=0. \quad (3.16)$$

This is *Seidel's identity* for the Genocchi numbers.

To obtain the expansion (3.2) we define $q-$Genocchi numbers $G_{2n}(q)$ by the generating function

$$z\frac{e(z)-e(-z)}{e(z)+e(-z)} = \sum_{n\geq 0}\frac{(-1)^{n-1}G_{2n}(q)(-q;q)_{2n-1}}{[2n]!}z^{2n}. \tag{3.17}$$

This implies that

$$t_{2n+1}(q) = \frac{G_{2n+2}(q)(-q;q)_{2n+1}}{[2n+2]}. \tag{3.18}$$

(Observe that this $q-$analogue of the Genocchi numbers does not coincide with the $q-$Genocchi numbers introduced by J. Zeng and J. Zhou which have been studied in [9]).

The first terms of the sequence $(G_{2n}(q))_{n\geq 1}$ are

$G_2(q) = 1,$

$G_4(q) = q\dfrac{1+q}{1+q^3},$

$G_6(q) = q^2\dfrac{(1+q)(1+q^2)(1+q+q^2)}{(1+q^4)(1+q^5)},$

$G_8(q) = q^3\dfrac{(1+q)^2(1+q^2)(1+q+3q^2+2q^3+3q^4+2q^5+3q^6+q^7+q^8)}{(1+q^5)(1+q^6)(1+q^7)}.$

**Theorem 3.3**

$$U_{2n+1}(x,s,q) = \sum_{k=0}^{n}\begin{bmatrix}2n+2\\2k\end{bmatrix}\frac{1}{[2k+1]}(-q;q)_{2n-2k+1}(-1)^{n-k}G_{2n-2k+2}(q)x^{2n+1-2k}U_{2k}(x,s,q). \tag{3.19}$$

**Proof**

In (2.38) we have seen that $U_n(1,s,q) = \sum_{k=0}^{\left\lfloor\frac{n}{2}\right\rfloor}\begin{bmatrix}n+1\\2k+1\end{bmatrix}q^{\binom{n-2k}{2}}(1+qs)(1+q^3s)\cdots(1+q^{2k-1}s).$

By comparing coefficients this is equivalent with

$$\frac{1}{e(-z)}\sum_{n\geq 0}\frac{z^{2n+1}}{[2n+1]!}(1+qs)(1+q^3s)\cdots(1+q^{2n-1}s) = \sum_{n\geq 1}\frac{U_{n-1}(1,s,q)}{[n]!}z^n. \tag{3.20}$$

Let now

$$U(z,s,q) = \sum_{n\geq 1}\frac{U_{n-1}(1,s,q)}{[n]!}z^n. \tag{3.21}$$

We then get

$$e(-z)U(z,s,q) = \sum_{n \geq 0} \frac{z^{2n+1}}{[2n+1]!}(1+qs)(1+q^3 s)\cdots(1+q^{2n-1}s) = -e(z)U(-z,s,q).$$

This implies
$$(e(z)-e(-z))(U(z,s,q)-U(-z,s,q)) = -e(z)U(-z,s,q) - e(-z)U(z,s,q) + e(z)U(z,s,q)$$
$$+e(-z)U(-z,s,q) = e(z)U(z,s,q) + e(-z)U(-z,s,q) = (e(z)+e(-z))(U(z,s,q)+U(-z,s,q)).$$

Since $U(z,s,q) + U(-z,s,q) = 2\sum_{n \geq 1} \frac{U_{2n-1}(1,s,q)}{[2n]!} z^{2n}$ and

$$U(z,s,q) - U(-z,s,q) = 2\sum_{n \geq 0} \frac{U_{2n}(1,s,q)}{[2n+1]!} z^{2n+1}$$

we see that

$$\frac{\sum_{n \geq 1} \frac{U_{2n-1}(1,s,q)}{[2n]!} z^{2n}}{\sum_{n \geq 0} \frac{U_{2n}(1,s,q)}{[2n+1]!} z^{2n+1}} = \frac{e(z)-e(-z)}{e(z)+e(-z)}. \tag{3.22}$$

Again the left-hand side does not depend on $s$. So we can e.g. choose $s = 0$ and get that

$$\frac{\sum_{n \geq 1} \frac{(-q;q)_{2n-1}}{[2n]!} z^{2n}}{\sum_{n \geq 0} \frac{(-q;q)_{2n}}{[2n+1]!} z^{2n+1}} = \frac{e(z)-e(-z)}{e(z)+e(-z)}. \tag{3.23}$$

If we write (3.22) in the form

$$\sum_{n \geq 1} \frac{U_{2n-1}(1,s,q)}{[2n]!} z^{2n} = z \frac{e(z)-e(-z)}{e(z)+e(-z)} \sum_{n \geq 0} \frac{U_{2n}(1,s,q)}{[2n+1]!} z^{2n}$$

and compare coefficients we get

$$U_{2n-1}(1,s,q) = \sum_{k=0}^{n} \begin{bmatrix} 2n \\ 2k \end{bmatrix} \frac{1}{[2k+1]} (-q;q)_{2n-2k-1} (-1)^{n-k-1} G_{2n-2k}(q) U_{2k}(1,s,q).$$

This immediately implies Theorem 3.3.

Since the left-hand side of (3.17) and $\frac{(-q;q)_{2n-1}}{[2n]!}$ are invariant under $q \to \frac{1}{q}$ we see that

$$G_{2n}\left(\frac{1}{q}\right) = G_{2n}(q). \tag{3.24}$$

Now we prove a $q$-analogue of (3.11):

**Theorem 3.4**

*The $q$ – Chebyshev polynomials $U_n(1,s,q)$ satisfy the identity*

$$\sum_{k=0}^{n}(-1)^k q^{\binom{k}{2}}\begin{bmatrix}n\\k\end{bmatrix}\prod_{j=n+m+1-k}^{n+m}(1+q^j)U_{2n+m-1-k}(1,s,q) = q^{n^2-n+mn}s^n U_{m-1}(1,s,q). \qquad (3.25)$$

**Proof**

Let

$$W(n,m,s,q) = \sum_{k=0}^{n}(-1)^k q^{\binom{k}{2}}\begin{bmatrix}n\\k\end{bmatrix}\prod_{j=n+m+1-k}^{n+m}(1+q^j)U_{2n+m-1-k}(1,s,q). \qquad (3.26)$$

We want to show that

$$W(n,m,s,q) == q^{n^2-n+mn}s^n U_{m-1}(1,s,q). \qquad (3.27)$$

We prove this identity with induction.
For $n=0$ it is the trivial identity $U_{m-1}(1,s,q) = U_{m-1}(1,s,q)$.
For $n=1$ it reduces to $U_{m+1}(1,s,q) - (1+q^{m+1})U_m(1,s,q) = q^m s U_{m-1}(1,s,q)$.
By definition of the polynomials this is true for all non-negative $m$.

In general we have

$$W(n,m,s,q) = W(n-1,m+2,s,q) - q^{n-1}(1+q^{m+1})W(n-1,m+1,s,q). \qquad (3.28)$$

Observing that

$$\begin{bmatrix}n-1\\k\end{bmatrix}(1+q^{n+m+1}) + q^{n-k}(1+q^{m+1})\begin{bmatrix}n-1\\k-1\end{bmatrix} =$$

$$= \left(\begin{bmatrix}n-1\\k\end{bmatrix} + q^{n-k}\begin{bmatrix}n-1\\k-1\end{bmatrix}\right) + q^{m+n+1-k}\left(\begin{bmatrix}n-1\\k-1\end{bmatrix} + q^k\begin{bmatrix}n-1\\k\end{bmatrix}\right) = \left(1+q^{m+n+1-k}\right)\begin{bmatrix}n\\k\end{bmatrix}$$

we get

$$W(n-1, m+2, s, q) - q^{n-1}(1+q^{m+1})W(n-1, m+1, s, q)$$

$$= \sum_{k=0}^{n-1}(-1)^k q^{\binom{k}{2}} \begin{bmatrix} n-1 \\ k \end{bmatrix} \prod_{j=n+m-k+2}^{n+m+1}(1+q^j)U_{2n+m-1-k}(1,s,q) - q^{n-1}(1+q^{m+1})\sum_{k=0}^{n-1}(-1)^k q^{\binom{k}{2}} \begin{bmatrix} n-1 \\ k \end{bmatrix} \prod_{j=n+m+1-k}^{n+m}(1+q^j)U_{2n+m-2-k}(1,s,q)$$

$$= \sum_{k=0}^{n-1}(-1)^k q^{\binom{k}{2}} \begin{bmatrix} n-1 \\ k \end{bmatrix} \prod_{j=n+m-k+2}^{n+m+1}(1+q^j)U_{2n+m-1-k}(1,s,q) - q^{n-1}(1+q^{m+1})\sum_{k=1}^{n}(-1)^{k-1} q^{\binom{k-1}{2}} \begin{bmatrix} n-1 \\ k-1 \end{bmatrix} \prod_{j=n+m+2-k}^{n+m}(1+q^j)U_{2n+m-1-k}(1,s,q)$$

$$= U_{2n+m-1}(1,s,q) + \sum_{k=1}^{n-1}(-1)^k q^{\binom{k}{2}} \prod_{j=n+m-k+2}^{n+m}(1+q^j)U_{2n+m-1-k}(1,s,q) \left( \begin{bmatrix} n-1 \\ k \end{bmatrix}(1+q^{n+m+1}) + q^{n-1-k+1}(1+q^{m+1}) \begin{bmatrix} n-1 \\ k-1 \end{bmatrix} \right)$$

$$+ q^{n-1}(1+q^{m+1})(-1)^n q^{\binom{n-1}{2}} \prod_{j=m+2}^{n+m}(1+q^j)U_{n+m-1}(1,s,q)$$

$$= \sum_{k=0}^{n}(-1)^k q^{\binom{k}{2}} \prod_{j=n+m-k+1}^{n+m}(1+q^j) \begin{bmatrix} n \\ k \end{bmatrix} U_{2n+m-1-k}(1,s,q) = W(n,m,s,q).$$

By induction (3.28) implies

$$W(n,m,s,q) = W(n-1, m+2, s, q) - q^{n-1}(1+q^{m+1})W(n-1, m+1, s, q)$$
$$= q^{n^2-n+(n-1)m} s^{n-1} U_{m+1}(1,s,q) - q^{n^2-n+(n-1)m} s^{n-1}(1+q^{m+1})U_m(1,s,q)$$
$$= q^{n^2-n+(n-1)m} s^{n-1} \left( U_{m+1}(1,s,q) - (1+q^{m+1})U_m(1,s,q) \right) = q^{n^2-n+nm} s^n U_{m-1}(1,s,q).$$

For $m = 0$ we get

$$\sum_{k=0}^{n}(-1)^k q^{\binom{k}{2}} \begin{bmatrix} n \\ k \end{bmatrix} \prod_{j=n-k+1}^{n}(1+q^j)U_{2n-1-k}(1,s,q) = 0. \tag{3.29}$$

An easy consequence is a $q$-analogue of the Seidel identity for the Genocchi numbers which gives an easy way to calculate the $q$-Genocchi numbers and shows that $(-q^{n+1}; q)_{n-1} G_{2n}(q) \in \mathbb{Z}[q]$ is a polynomial with integer coefficients.

**Theorem 3.5 (q-Seidel formula)**

$$\sum_{k=0}^{\lfloor \frac{n}{2} \rfloor} q^{\binom{2k}{2}} \begin{bmatrix} n \\ 2k \end{bmatrix} (-1)^k \frac{(-q^{n-2k+1}; q)_{2k}}{(-q^{2n-2k}; q)_{2k}} G_{2n-2k}(q) = [n = 1]. \tag{3.30}$$

**Proof**

Since the set of polynomials $\{U_{2n}(1,s,q)\}_{n \geq 0}$ is a basis for the vector space of polynomials in $s$ we can define a linear functional $\lambda$ by

$$\lambda(U_{2n}(1,s,q)) = [n = 0]. \tag{3.31}$$

By (3.19) this implies

$$\lambda\left(U_{2n-1}(1,s,q)\right) = (-1)^{n-1}\left(-q;q\right)_{2n-1} G_{2n}(q). \tag{3.32}$$

If we apply this to (3.29) we get for $n>1$

$$0 = \lambda\left(\sum_{k=0}^{n}(-1)^k q^{\binom{k}{2}} \begin{bmatrix} n \\ k \end{bmatrix} \prod_{j=n-k+1}^{n}(1+q^j) U_{2n-1-k}(1,s,q)\right) = \sum_{k=0}^{n}(-1)^k q^{\binom{k}{2}} \begin{bmatrix} n \\ k \end{bmatrix} \prod_{j=n-k+1}^{n}(1+q^j) \lambda(U_{2n-1-k}(1,s,q))$$

$$\sum_{k=0}^{\lfloor \frac{n}{2} \rfloor} q^{\binom{2k}{2}} \begin{bmatrix} n \\ 2k \end{bmatrix} \prod_{j=n-k+1}^{n}(1+q^j) \prod_{j=1}^{2n-1-2k}(1+q^j)(-1)^{n-k-1} G_{2n-2k}(q).$$

Dividing by $(-q;q)_{2n-1}$ we get (3.30).

It should be noted that just as for $q=1$ (3.30) is in fact the same formula as (3.14). We need only use (3.18) to translate one formulation into the other.

Finally we want to show how to derive a Seidel triangle for the $q-$Genocchi numbers. We construct the following triangle consisting of numbers $a(n,k,q)$ with $n=0,1,2,\cdots$ and $0 \le k \le 1 + \left\lfloor \dfrac{n}{2} \right\rfloor$.

Let $a(2n,k,q) = (-1)^n s^{n+1-k} U_{2k-2}(1,s,q)$ and $a(2n+1,k,q) = (-1)^n s^{n+1-k} U_{2k-1}(1,s,q)$.

The first terms are (if we delete the column $k=0$)

$$\begin{array}{llll}
U_0(1,s,q) & & & \\
U_1(1,s,q) & & & \\
-sU_0(1,s,q) & -U_2(1,s,q) & & \\
-sU_1(1,s,q) & -U_3(1,s,q) & & \\
s^2 U_0(1,s,q) & sU_2(1,s,q) & U_4(1,s,q) & \\
s^2 U_1(1,s,q) & sU_3(1,s,q) & U_5(1,s,q) & \\
-s^3 U_0(1,s,q) & -s^2 U_2(1,s,q) & -sU_4(1,s,q) & U_6(1,s,q) \\
\vdots & \vdots & \vdots & \vdots
\end{array}$$

Then
$$a(2n+1,k,q) = q^{2k-2} a(2n+1,k-1,q) + (1+q^{2k-1}) a(2n,k,q)$$
for $k = 1, 2, \cdots, n+1$.

On the other hand
$$a(2n,k,q) = q^{1-2k}\left(a(2n,k+1,q) + (1+q^{2k}) a(2n-1,k,q)\right)$$
for $k = 1, 2, \cdots, n$.
For $k = n+1$ we get $a(2n, n+1, q) = U_{2n}(1,s,q)$.

If we apply the linear functional $\lambda$ and let $b(n,k,q) = \lambda(a(n,k,q))$ then $b(2n,n+1,q) = 0$ and therefore we have $b(2n,n+1,q) = q^{1-2k}\left(b(2n,n+2,q) + (1+q^{2k})b(2n-1,n+1,q)\right) = 0$.

Thus we get

**Theorem 3.6 (q-Genocchi triangle)**

*Define a triangle $(b(n,k,q))$ for $n \in \mathbb{N}$ and $0 \leq k \leq 1 + \left\lfloor \dfrac{n}{2} \right\rfloor$ by*

$$b(2n+1,k,q) = q^{2k-2}b(2n+1,k-1,q) + (1+q^{2k-1})b(2n,k,q) \tag{3.33}$$

*and*

$$b(2n,k,q) = q^{1-2k}\left(b(2n,k+1,q) + (1+q^{2k})b(2n-1,k,q)\right) \tag{3.34}$$

*for $1 \leq k \leq n+1$ with initial values $b(0,1,q) = 1$ and $b(1,1,q) = 1+q$.*
*Then*

$$b(2n-1,n) = \lambda\left((-1)^{n-1}U_{2n-1}(1,s,q)\right) = (-q;q)_{2n-1} G_{2n}(q). \tag{3.35}$$

This is another simple method to compute the $q-$Genocchi numbers.